\documentclass[12pt,a4paper,oneside,reqno]{article}

\usepackage[utf8]{inputenc}
\usepackage[english]{babel}
\usepackage[left=1in,right=1in,top=1in,bottom=1in]{geometry}
\usepackage{amsfonts,amssymb,amsmath,amsthm}
\usepackage[symbol]{footmisc}
\usepackage{longtable,array}
\newcolumntype{L}[1]{>{\raggedright\arraybackslash$}p{#1}<{$}}
\usepackage{enumerate}
\usepackage{cite}
\usepackage{xcolor}

\usepackage{hyperref}
\hypersetup{hidelinks}
\usepackage{fancyhdr}
\theoremstyle{plain}
\newtheorem{theorem}{Theorem}

\theoremstyle{definition}

\makeatletter
\@removefromreset{theorem}{section}\@removefromreset{equation}{section}

\makeatother
\usepackage{adjustbox}
\usepackage{tikz}\usetikzlibrary{arrows}

\IfFileExists{fouriernc.sty}{\usepackage{fouriernc}}{}

\newcommand{\doi}[1]{\href{https://doi.org/#1}{doi:\detokenize{#1}}}
\newcommand{\arxiv}[1]{\href{https://arxiv.org/abs/#1}{arXiv:#1}}

\theoremstyle{plain}
\newtheorem{Th}[theorem]{Theorem}
\newtheorem{Lem}[theorem]{Lemma}
\theoremstyle{definition}

\begin{document}

\noindent\fbox{\parbox{\dimexpr\textwidth-2\fboxsep-2\fboxrule\relax}{\small
\textbf{Corrected version (verified by Claude).} This is a corrected version of the paper I.~Kaygorodov, \emph{On the Kantor product}, Journal of Algebra and its Applications, 16 (2017), no.~9, 1750167, \doi{10.1142/S0219498817501675}.}}

\bigskip

\begin{center}
{\Large On the Kantor product}

\bigskip

{\bf Ivan Kaygorodov}
\end{center}

\medskip

\noindent{\small\textbf{Abstract.} We study the algebra of bilinear multiplications of an $n$-dimensional vector space. In particular, we study the Kantor product of some well-known multiplications (associative, Lie, alternative, Novikov, and some others).}

\bigskip

\noindent {\bf Keywords:}
{\it Kantor product, conservative algebra, bilinear multiplication, polynomial identity, alternative algebra.}

\bigskip

\noindent {\bf MSC2020:} 17A30 (primary); 16R10, 17A32, 17A36, 17D05 (secondary).

\bigskip

\tableofcontents

\bigskip

\section*{Introduction}\addcontentsline{toc}{section}{Introduction}

Kantor introduced the class of conservative algebras in \cite{Kantor72}.
This class includes some well-known algebras: associative, Jordan, Lie, Leibniz and Zinbiel
\cite{KLP}.
The conservative algebra $U(n)$ is of great importance in the theory of conservative algebras~\cite{Kantor90}:
it plays a role analogous to the role of $\mathfrak{gl}_n$ in the theory of Lie algebras.
The space of the algebra $U(n)$ is the space of all bilinear multiplications on
the $n$-dimensional space $V_n$.
To define  the operation of multiplication $[ \ , \ ]$ in the algebra $U(n)$
we fix a vector $u \in V_n$ and
for two multiplications $A, B \in U(n)$ and two elements $x,y \in V_n$ we  set
\begin{eqnarray}\label{Kanpro}
x*y=[A,B](x,y) = A(u,B(x,y))-B(A(u,x),y)-B(x,A(u,y)).
\end{eqnarray}
Some properties of the algebra $U(2)$ were studied in \cite{KLP,KV15}.
We say that the product of two multiplications on $n$-dimensional vector space
defined by (\ref{Kanpro}), is the left Kantor product of these multiplications.
In a similar way we can define the right Kantor product and obtain similar results.
We assume that the Kantor product is the left Kantor product.
The Kantor product of a multiplication  $\cdot$ by itself is {\it the Kantor square} of $\cdot$ and it is denoted by
 $[\cdot, \cdot].$
It gives us a map $K$ from any variety $V$ of algebras to some class $K(V).$

The Kantor square of a multiplication $\cdot$ can be rewritten (see \cite{Kantor72})
 as the product of the left multiplication $L_u$
and the multiplication $\cdot,$ as $[L_u, \cdot],$ where
$$[L_u, \cdot](x,y) = u\cdot(x\cdot y)- (u\cdot x)\cdot y - x \cdot (u\cdot y)=[\cdot, \cdot](x,y).$$
The multiplication $[L_u, \cdot]$ plays an important role in the definition of (left) conservative algebras
 \cite{kacan,Kantor72}.
We recall that an algebra $A$ with a multiplication $\cdot$ is called a (left) conservative algebra
if and only if there exists a new multiplication $*$ such that
$$[L_a,[L_b, \cdot]]=- [L_{a*b},\cdot].$$

The main aim of this paper is to study the properties of the Kantor product of multiplications.
One of the central questions studied in this paper
is the following:

{\it\bf Question.} What identities does the class of algebras $K(V)$ satisfy if we know the identities of $V$?

We give some partial answers to this question for
associative, (anti)commutative, Perm, Lie, Leibniz, Zinbiel, left-commutative, bicommutative,
Novikov, alternative, quasi-associative and quasi-alternative algebras;
 we also describe the Kantor product of multiplications
in associative dialgebras, duplicial, dual duplicial, $As^{(2)}$,
Poisson, generalized Poisson and Novikov-Poisson algebras.
Finally,  we study the Kantor square in some special cases;
in particular,
the associative algebras with identities, nilpotent and right-nilpotent algebras,
associative algebras isomorphic to their Kantor squares;
and discuss the coincidence of derivations and automorphisms of the algebra and its Kantor square.
We formulate the following

{\it\bf Open problem.} Is $K(V)$ a variety of algebras for some variety $V$?

    \section{The Kantor square}

\

In this section we omit the technical and straightforward proofs of the lemmas.
We are using the standard notation:
$$(a,b,c)_*=(a*b)*c-a*(b*c), (a,b,c)=(ab)c-a(bc); $$
$$[a,b]_*=a*b-b*a, [a,b]=ab-ba;$$
$$\circlearrowright_{a,b,c} [f(a,b,c)]=f(a,b,c)+f(b,c,a)+f(c,a,b).$$

\subsection{Associative algebras.}
The variety of {\it associative algebras} is defined by the identity
$$(xy)z=x(yz).$$

    \begin{Lem}\label{thass} Let $(A;\cdot)$ be an associative algebra.
    Then $(A;[\cdot, \cdot])$ is an associative algebra.
    \end{Lem}

\

\subsection{(Anti)commutative algebras.}
The variety of {\it  $($anti$)$commutative algebras} is defined by the identity
$$xy=\epsilon yx,$$ where $\epsilon=1$ in the commutative case
 and $\epsilon =-1$ in the  anticommutative case.

    \begin{Lem}\label{anticom} Let $(A;\cdot)$ be an $($anti$)$commutative algebra.
    Then $(A;[\cdot, \cdot])$ is an $($anti$)$commutative algebra.
    \end{Lem}

\

\subsection{Perm algebras.}
The variety of {\it Perm  algebras} (see, for example, \cite{Chapa01}) is defined by the identity
$$(xy)z=x(yz)=x(zy).$$

    \begin{Lem} Let $(A;\cdot)$ be a Perm algebra.
    Then $(A;[\cdot, \cdot])$ is a Perm  algebra.
    \end{Lem}

\

\subsection{Lie algebras.}
The variety of {\it Lie algebras} is defined by the identities
$$xy=-yx, (xy)z+(yz)x+(zx)y=0.$$

    \begin{Lem}Let $(A;\cdot)$ be a Lie algebra.
    Then $[\cdot, \cdot]=0.$
    \end{Lem}

\

\subsection{Leibniz algebras.}
The variety of (left) {\it Leibniz algebras} (see, for example, \cite{Dzhuma01})
 is defined by the identity
$$x(yz)=(xy)z+y(xz).$$

    \begin{Lem}\label{leibniz}Let $(A;\cdot)$ be a $($left$)$ Leibniz algebra.
    Then $[\cdot, \cdot]=0.$
    \end{Lem}

\

\subsection{Left-commutative algebras.}
The variety of {\it  left-commutative algebras} (see, for example, \cite{shestakov_kuzm})
 includes commutative-associative, bicommutative, Novikov, Zinbiel algebras and some other.
This variety is defined by the identity
$$x(yz)=y(xz).$$

    \begin{Lem}\label{thleftcom} Let $(A;\cdot)$ be a left-commutative algebra.
    Then $(A;[\cdot, \cdot])$ is a left-commutative algebra.
    \end{Lem}

\

\subsection{Bicommutative algebras.}
The variety of {\it bicommutative algebras} (see, for example, \cite{Dzhuma11})
 is defined by the identities
$$x(yz)=y(xz), (xy)z=(xz)y.$$

    \begin{Lem} Let $(A;\cdot)$ be a bicommutative algebra.
    Then $(A;[\cdot, \cdot])$ is an associative-commutative algebra.
    \end{Lem}

\

\subsection{Zinbiel algebras.}
The variety of (left) {\it  Zinbiel algebras} (see, for example, \cite{Dzhuma07})
is defined by the identity
$$x(yz)=(xy)z+(yx)z.$$

    \begin{Lem} Let $(A;\cdot)$ be a $($left$)$ Zinbiel algebra.
    Then $(A;[\cdot, \cdot])$ is a $($left$)$ Zinbiel algebra.
    \end{Lem}

\

\subsection{Novikov algebras.}
The variety of (left) {\it Novikov algebras} (see, for example, \cite{Dzhuma14})
 is defined by the identities
$$x(yz)=y(xz), (x,y,z)=(x,z,y).$$

    \begin{Lem}Let $(A;\cdot)$  be a $($left$)$ Novikov algebra.
    Then  $(A;[\cdot, \cdot])$ is a $($left$)$ Novikov algebra.
    \end{Lem}

\

\subsection{Alternative algebras.}
The variety of {\it alternative algebras} (see, for example, \cite{KP})
is defined by the identities
\begin{eqnarray}\label{altern}
x^2y=x(xy),\qquad xy^2=(xy)y.
\end{eqnarray}
It is also well known (see, for example, \cite{rings})
that an  alternative algebra is flexible:  $$(xy)x=x(yx);$$
and satisfies the Moufang identities:
$$x(yzy)=((xy)z)y, (yzy)x=y(z(yx)), (xy)(zx)=x(yz);$$
and the following identities hold:
$$(x,y,z)=-(y,x,z), (x,y,z)=-(x,z,y),$$
$$(x,xy,z)=(x,y,z)x,(x,yx,z)=x(x,y,z).$$

The main example of a non-associative alternative algebra is
 a Cayley --- Dickson algebra ${\bf C}$ \cite{rings}.
Let $F$ be a field of characteristic $\neq 2.$
It is an algebra ${\bf C}$ with the basis $e_0 = 1, e_1,\dots,e_7$
and the following multiplication table:
\begin{equation}
\label{table}
\begin{array}{|c||c|c|c|c|c|c|c|} \hline
1 &  e_1 & e_2 & e_3 & e_4 & e_5 & e_6 & e_7 \\
\hline \hline e_1 & \alpha \cdot 1 &  e_3  &  \alpha e_2 & e_5 &  \alpha e_4  & -e_7 &  -  \alpha e_6 \\
\hline
 e_2  & -e_3 & \beta \cdot 1 & - \beta e_1 & e_6 & e_7 &  \beta  e_4 &  \beta e_5 \\
\hline
 e_3  & - \alpha e_2 &  \beta e_1 & -\alpha \beta \cdot 1 & e_7 &  \alpha e_6 &  -\beta e_5 &
 -\alpha \beta  e_4 \\
\hline
 e_4 & -e_5 & -e_6 & -e_7 & \gamma \cdot 1 & - \gamma e_1 & -\gamma e_2 &  -\gamma e_3 \\
\hline
 e_5 & - \alpha e_4 & -e_7 & -\alpha  e_6 & \gamma e_1 & -\alpha\gamma \cdot 1  &  \gamma  e_3 &
  \alpha \gamma e_2 \\
\hline
 e_6 & e_7 & -\beta e_4 &  \beta e_5 & \gamma e_2 & -\gamma e_3 &  -\beta \gamma \cdot 1 &
 -\beta\gamma  e_1 \\
\hline
 e_7 & \alpha  e_6 & -\beta e_5 &  \alpha \beta e_4 &  \gamma e_3 &  -\alpha \gamma e_2 &
 \beta \gamma e_1 & \alpha\beta\gamma \cdot 1 \\
\hline
\end{array}
\end{equation}

    \begin{Th}\label{teoremaalt}Let $(A;\cdot)$  be an alternative algebra.
    Then  $(A;[\cdot, \cdot])$ is a flexible algebra. Furthermore,\\
    1$)$ $(A;[\cdot, \cdot])$ is an alternative algebra if and only if $A$ satisfies the identity
    \begin{eqnarray}\label{altass}(x,u,(x,u,y))=0;\end{eqnarray}
    2$)$ $(A;[\cdot, \cdot])$ is a noncommutative Jordan algebra if and only if $A$
    satisfies the identity
    $$((xuxu)y)(ux)-(xuxu)(y(ux))+((x(u(x(uy))))u)x-x(u(x(u((yu)x))))=0$$
    (here $xuxu=(xu)(xu)$ is well defined, since the subalgebra generated by $x$ and $u$ is associative by Artin's theorem)$;$\\
    3$)$ $(A;[\cdot, \cdot])$ is a Jordan algebra if $(A, \cdot)$ is a commutative alternative algebra$;$\\
    4$)$ $({\bf C};[\cdot, \cdot])$ is alternative for a Cayley --- Dickson algebra ${\bf C}$,
    if and only if $u= u_0\cdot 1.$
    \end{Th}

{\bf Proof.}
It is easy to see that
$$a * b = u(ab)-(ua)b-a(ub)= (au)b- 2a(ub).$$

Now, we can see
$$ (x*y)*x-x*(y*x)= ((xu)y-2x(uy))*x- x*((yu)x-2y(ux))=$$
$$(((xu)y)u)x-2((x(uy))u)x-2((xu)y)(ux)+4(x(uy))(ux)-$$
$$(xu)((yu)x)+2(xu)(y(ux))+2x(u((yu)x))-4x(u(y(ux)))=$$
$$x(uyu)x-2((x,u,y)u)x-2x(uyu)x-2x(uyu)x+2((xu)y,u,x)+4x(uyu)x-$$
$$x(uyu)x-2(xy)(y,u,x)+2x(uyu)x+2x(uyu)x+2x(u(y,u,x))-4x(uyu)x=$$
$$2[((xu)y,u,x)-(x,u,x(uy))-((xy)(y,u,x)-x(u(y,u,x)))]=$$
$$2[((x,u,y),u,x)-((x,u,y),u,x)]=0.$$

It follows that $(A,*)$ is a flexible algebra.

$1)$ It is easy to see that a flexible algebra is alternative if and only if it satisfies the first identity from (\ref{altern}).
We have
$$(x*x)*y - x*(x*y)=$$
$$-(xuxu)y+2(xux)(uy)
-(xu)((xu)y)+2x(u((xu)y))+2(xu)(x(uy))-4x(u(x(uy)))=$$
$$-2(xu)^2y-2x(u(x(uy)))+2x(u((xu)y))+2(xu)(x(uy))=-2(x,u,(x,u,y)).$$
Now, since the characteristic of the ground field is different from $2$, the multiplication $[\cdot, \cdot]$ is alternative if and only if $(x,u,(x,u,y))=0.$

$2)$ Note that an algebra $B$ is a {\it  non-commutative Jordan algebra}
if and only if it is flexible and it satisfies the Jordan identity:
$(x^2,y,x)=0.$

Obviously,
$$((x*x)*y)*x-(x*x)*(y*x)=$$
$$(-(xuxu)y+2(xux)(uy))*x + (xux)*((yu)x-2y(ux))=$$
$$-(((xuxu)y)u)x+2((xuxu)y)(ux)+2(((xux)(uy))u)x-4((xux)(uy))(ux)+$$
$$(xuxu)((yu)x)-2(xux)(u((yu)x))-2(xuxu)(y(ux))+4(xux)(u(y(ux)))=$$
$$2(((xuxu)y)(ux)-(xuxu)(y(ux)) +((x(u(x(uy))))u)x-x(u(x(u((yu)x)))).$$

Now, since $(A,*)$ is flexible, $[\cdot, \cdot]$ is a noncommutative Jordan multiplication if and only if
$((xuxu)y)(ux)-(xuxu)(y(ux))+((x(u(x(uy))))u)x-x(u(x(u((yu)x))))=0.$

$3)$ It is easy to see that
if $A$ is a commutative alternative algebra then we have
$$((xuxu)y)(ux)-(xuxu)(y(ux)) +(((xux)(uy))u)x-(xux)(u((yu)x))=$$
$$(xu)((xu)((xu)y)))- (xu)((xu)((xu)y)))+
x(u(x(u(x(uy)))))-x(u(x(u(x(uy)))))=0.$$
It follows that $[\cdot, \cdot]$ is  non-commutative Jordan
and from Lemma~\ref{anticom} we infer that $[\cdot, \cdot]$ is  Jordan.

$4)$ By $1)$, $({\bf C},*)$ is alternative if and only if $(x,u,(x,u,y))=0$ for all $x,y\in{\bf C}$. If $u=u_0\cdot 1$, then this identity obviously holds. Conversely, let $u=u_0\cdot 1+u_1e_1+\ldots+u_7e_7$ satisfy~(\ref{altass}). Since an associator with an argument from $F\cdot 1$ is zero, we may assume that $u_0=0$. Let $c_k=e_k^2\in F$, that is,
$$c_1=\alpha,\ c_2=\beta,\ c_3=-\alpha\beta,\ c_4=\gamma,\ c_5=-\alpha\gamma,\ c_6=-\beta\gamma,\ c_7=\alpha\beta\gamma.$$
For distinct $i,j\in\{1,\ldots,7\}$ let $m$ be the index such that $e_ie_j\in F e_m$. A direct computation with the multiplication table shows that the coefficient of $e_j$ in $(e_i,u,(e_i,u,e_j))$ is equal to
$$-4c_i\sum_{k\notin\{i,j,m\}}c_ku_k^2.$$
The triples $\{i,j,m\}$ are exactly the seven triples $\{1,2,3\}$, $\{1,4,5\}$, $\{1,6,7\}$, $\{2,4,6\}$, $\{2,5,7\}$, $\{3,4,7\}$, $\{3,5,6\}$ (the lines of the Fano plane). Since $c_i\neq0$ and the characteristic is different from $2$, we obtain the following system:
\begin{equation*}
\begin{array}{ccccccccc}
\alpha u^2_1 &              & -\alpha\beta u^2_3 & + \gamma u^2_4 &                      &
-\beta\gamma u^2_6 &                          & =& 0, \\
\alpha u^2_1 & +\beta u^2_2 &                    &                & -\alpha \gamma u_5^2  &
-\beta\gamma u^2_6 &                          & =& 0, \\
\alpha u^2_1 & +\beta u^2_2 &                    & + \gamma u^2_4 &                      &
& +\alpha\beta\gamma u_7^2  & =& 0, \\
\alpha u^2_1 &              & -\alpha\beta u^2_3 &                & -\alpha \gamma u_5^2  &
                  & +\alpha\beta\gamma u_7^2  & =& 0, \\
             & +\beta u^2_2 & -\alpha\beta u^2_3 & + \gamma u^2_4 & -\alpha \gamma u_5^2  &
                               &                          & =& 0, \\
             & +\beta u^2_2 & -\alpha\beta u^2_3 &                &                      &
             -\beta\gamma u^2_6 & +\alpha\beta\gamma u_7^2  & =& 0, \\
             &              &                    & + \gamma u^2_4 & -\alpha \gamma u_5^2  &
             -\beta\gamma u^2_6 & +\alpha\beta\gamma u_7^2  & =& 0.  \\
\end{array}
\end{equation*}
This is a linear system in the unknowns $w_k=c_ku_k^2$ whose matrix $M$ has the entry $1$ in the position $(\ell,k)$ if $k\notin\ell$ and $0$ otherwise. We have $M=J-N$, where $J$ is the all-ones matrix and $N$ is the incidence matrix of the Fano plane. Since $NN^T=2I+J$, the matrix $M$ acts as the multiplication by $4$ on the all-ones vector and as $-N$ on its orthogonal complement, and $\det M=-32\neq0$. Hence $w_k=c_ku_k^2=0$, so $u_1=\ldots=u_7=0$ and $u= u_0\cdot 1.$

The theorem is proved.

\subsection{Quasi-associative algebras.}
{\it Quasi-associative algebras} (see, for example, \cite{dedkov89})
are defined by the identities
\begin{center}$(x,y,z)+(y,z,x)+(z,x,y)=0,$\qquad 
$(x,y,z) = \alpha [y,[x,z]],$\end{center}
where $\alpha$ is a fixed element in the ground field $F.$
It is known \cite{dedkov89} that an algebra $(A, \cdot)$ is quasi-associative
if and only if there exists an associative algebra
$A$ with a new multiplication such that for some $\lambda\in F$:
$$x\cdot y = \lambda xy + (1-\lambda) yx.$$

    \begin{Lem} Let $(A;\cdot)$ be a quasi-associative algebra.
    Then $(A;[\cdot, \cdot])$ is a quasi-associative algebra.
    \end{Lem}

\

\subsection{Quasi-alternative algebras.}
{\it Quasi-alternative algebras} (see, for example, \cite{dedkov89})
are defined by the identities
$$(x,y,x)=0,$$
$$(x,x,y) = \alpha [x,[x,y]],$$
where $\alpha$ is a fixed element from the ground field $F.$
It is known \cite{dedkov89} that an algebra $(A, \cdot)$ is a quasi-alternative algebra
if and only if there exists an alternative algebra $A$ with a new multiplication such that for some $\lambda\in F$:
$$x\cdot y = \lambda xy + (1-\lambda) yx.$$

    \begin{Lem} Let $(A;\cdot)$ be a quasi-alternative algebra.
    Then $(A;[\cdot, \cdot])$ is a flexible algebra.
    \end{Lem}

\

\subsection{Associative dialgebras.}
The variety of {\it associative dialgebras} (see, for example, \cite{kolesn08})
is defined by the identities
$$(x \vdash y) \vdash z = ( x \dashv y)  \vdash z, x \dashv ( y \vdash z) =  x \dashv (y   \dashv z),$$
$$(x \vdash y) \vdash z = x \vdash (y \vdash z), (x \dashv y) \dashv z = x \dashv (y \dashv z),
(x \vdash y) \dashv z = x \vdash (y \dashv z).$$

    \begin{Lem}Let $(A;\vdash,\dashv)$ be an associative dialgebra.
    Then  $(A;[\vdash, \dashv])$  is an associative algebra.
    \end{Lem}

\

\subsection{Duplicial algebras.}

The variety of {\it duplicial algebras} (see, for example, \cite{Loday08}) is defined by the identities
\begin{center}
    $(x \prec y) \prec z = x \prec (y \prec z ),$ \quad
$(x \succ y) \prec z = x \succ (y \prec z),$ \quad 
$(x \succ y) \succ z = x \succ (y \succ z).$
\end{center}

    \begin{Lem}Let $(A;\prec,\succ)$ be a duplicial algebra.
    Then  $(A;[\succ ,\prec ])$  is an associative algebra.
    \end{Lem}

\subsection{Dual duplicial algebras.}

The variety of {\it dual duplicial algebras} (see, for example, \cite{zinbiel})
 is defined by the identities
$$(x \prec y) \prec z = x \prec (y \prec z ),(x \succ y) \prec z = x \succ (y \prec z),(x \succ y)
\succ z = x \succ (y \succ z),$$
$$ x \prec  (y \succ z) = ( x \prec  y) \succ z =0.$$

    \begin{Lem}Let $(A;\prec,\succ)$ be a dual duplicial algebra.
    Then  $[\succ ,\prec ]=0$ and $(A;[\prec , \succ ])$  is a $2$-nilpotent algebra.
    \end{Lem}

\subsection{$As^{(2)}$-algebras.}

The variety of {\it  $As^{(2)}$-algebras} (see, for example, \cite{zinbiel}) is defined by the identities
$$(x \circ y) \cdot z = x \circ (y \cdot z), (x \cdot y) \circ  z = x \cdot (y \circ z),  $$
$$(x \circ y) \circ z = x \circ (y \circ z), (x \cdot y) \cdot  z = x \cdot (y \cdot z).  $$

    \begin{Lem}Let $(A;\cdot,\circ)$ be a $As^{(2)}$-algebra.
    Then  $(A;[\cdot,\circ ])$ and $(A;[\circ  , \cdot ])$  are associative algebras.
    \end{Lem}

\subsection{Commutative tridendriform algebra.}
The variety of {\it commutative tridendriform algebras} (see, for example, \cite{Loday07})
is defined by the identities
$$x \cdot y  = y \cdot x, (x \cdot y) \cdot z = x \cdot (y \cdot z),$$
$$(x \prec y) \prec z = x \prec (y \prec z ) + x \prec (z \prec y ),$$
$$(x\cdot y)  \prec z=x \cdot (y  \prec z).$$

    \begin{Lem}Let $(A;\cdot,\prec)$ be a commutative tridendriform algebra.
    Then   $(A;[\prec,\cdot ])$ is a commutative algebra and  $(A;[\cdot  , \prec ])$
    is a right Zinbiel algebra.
    \end{Lem}

\subsection{Poisson algebras.}
The variety of {\it Poisson algebras} (see, for example, \cite{shestakov93}) is defined by the identities
$$xy=yx, (xy)z=x(yz), \{xy,z\}=\{x,z\}y+x\{y,z\},$$
$$\{x,y\}=-\{y,x\}, \{\{x,y\},z\}+\{\{y,z\},x\}+\{\{z,x\},y\}=0.$$

    \begin{Th}\label{ThPoisson} Let $(A;\cdot,\{,\})$ be a Poisson algebra.
    Then  $[\{,\}, \cdot]=0$ and  $(A;[\cdot, \{,\}])$ is a Lie  algebra.
    \end{Th}

{\bf Proof.}
In the first case, we have
$$a*b= \{u,ab\} - \{u,a\}b-a\{u,b\}=0.$$

In the second case,
$$a*b= u\{a,b\}-\{ua,b\}-\{a,ub\}= - a\{u,b\}- b\{a,u\}-u\{a,b\}=-b*a.$$
and
$$(a*b)*c+(b*c)*a+(c*a)*b=$$
$$-(a\{u,b\}+\{a,u\}b+\{a,b\}u)*c$$
$$  -(b\{u,c\}+\{b,u\}c+\{b,c\}u)*a$$
$$  -(c\{u,a\}+\{c,u\} a+\{c,a\}u)*b=$$

$$  a\{u,b\}\{u,c\}+\{a\{u,b\},u\}c+\{a\{u,b\},c\}u+$$
$$    \{a,u\}b\{u,c\}+\{\{a,u\}b,u\}c+\{\{a,u\}b,c\}u+$$
$$ \{a,b\}u\{u,c\}+\{\{a,b\}u,u\}c+\{\{a,b\}u,c\}u+$$
$$    b\{u,c\}\{u,a\}+\{b\{u,c\},u\}a+\{b\{u,c\},a\}u+$$
$$ \{b,u\}c\{u,a\}+\{\{b,u\}c,u\}a+\{\{b,u\}c,a\}u+$$
$$    \{b,c\}u\{u,a\}+\{\{b,c\}u,u\}a+\{\{b,c\}u,a\}u+$$
$$ c\{u,a\}\{u,b\}+\{c\{u,a\},u\}b+\{c\{u,a\},b\}u+$$
$$   \{c,u\}a\{u,b\}+\{\{c,u\}a,u\}b+\{\{c,u\}a,b\}u+$$
$$\{c,a\}u\{u,b\}+\{\{c,a\}u,u\}b+\{\{c,a\}u,b\}u=$$

$$a\{ u,b\}\{ u,c\}+\{ u,b\}\{ a,u\}c+ca\{\{ u,b\} ,u\}+au\{ \{ u,b\},c\}+\{ u,b\}\{ a,c\}u+$$
$$\{ a,u\}b\{ u,c\}+\{ \{ a,u\},u\}bc+\{a ,u\}\{ b,u\}c+bu\{\{a ,u\} ,c\}+\{ a,u\}\{ b,c\}u+$$
$$\{ a,b\}u\{ u,c\}+\{ \{ a,b\},u\}uc+\{ a,b\}\{ u,c\}u+uu\{ \{ a,b\},c\}+$$
$$b\{ u,c\}\{ u,a\}+\{ b,u\}\{ u,c\}a+ba\{ \{ u,c\},u\}+\{ \{ u,c\},a\}bu+\{ u,c\}\{ b,a\}u+$$
$$\{ b,u\}c\{ u,a\}+\{ b,u\}\{ c,u\}a+ca\{ \{ b,u\},u\}+\{ b,u\}\{ c,a\}u+\{ \{ b,u\},a\}cu+$$
$$\{ b,c\}u\{ a,u\}+\{ \{ b,c\},u\}ua+\{ b,c\}\{ u,a\}u+uu\{ \{ b,c\},a\}+$$
$$c\{ u,a\}\{ u,b\}+\{ c,u\}\{ u,a\}b+cb\{ \{ u,a\},u\}+\{ c,b\}\{ u,a\}u+cu\{ \{ u,a\},b\}+$$
$$\{ c,u\}a\{ u,b\}+\{ \{ c,u\},u\}ab+\{ c,u\}\{ a,u\}b+\{ c,u\}\{ a,b\}u+\{ \{ c,u\},b\}au+$$
$$\{ c,a\}u\{ u,b\}+\{ \{ c,a\},u\}ub+\{ c,a\}\{ u,b\}u+uu\{ \{ c,a\},b\}=$$

$$  (a\{ u,b\}\{ u,c\}+ \{ c,u\}a\{ u,b\})+
    (\{ u,b\}\{ a,u\}c+ \{a ,u\}\{ b,u\}c)+$$
$$  (ca\{\{ u,b\} ,u\}+ ca\{ \{ b,u\},u\})+
    (\{ u,b\}\{ a,c\}u+ \{ c,a\}\{ u,b\}u)+$$
$$  (\{ a,u\}b\{ u,c\}+ b\{ u,c\}\{ u,a\})+
    (\{ \{ a,u\},u\}bc+ cb\{ \{ u,a\},u\})+$$
$$  (\{ a,u\}\{ b,c\}u+ \{ b,c\}u\{ a,u\})+
    (\{ a,b\}u\{ u,c\}+ \{ c,u\}\{ a,b\}u)+$$
$$  (\{ a,b\}\{ u,c\}u+ \{ b,u\}\{ u,c\}a)+
    (ba\{ \{ u,c\},u\}+ \{ \{ c,u\},u\}ab)+$$
$$  (\{ b,u\}c\{ u,a\}+ \{ c,u\}\{ u,a\}b)+
    (\{ b,u\}\{ c,a\}u+ \{ c,a\}u\{ u,b\})+$$
$$  (\{ b,c\}\{ u,a\}u+ \{ c,b\}\{ u,a\}u)+
    (\{ c,u\}\{ u,a\}b+ \{ c,u\}\{ a,u\}b)+$$
$$  [au\{ \{ u,b\},c\}+au\{ \{ b,c\},u\}+au\{ \{ c,u\},b\}]+$$
$$  [bu\{\{a ,u\} ,c\}+bu\{ \{ u,c\},a\}+bu\{ \{ c,a\},u\}]+$$
$$  [cu\{ \{ a,b\},u\}+cu\{ \{ u,a\},b\}+cu\{ \{ b,u\},a\}]+$$
$$  [uu\{ \{ a,b\},c\}+uu\{ \{ b,c\},a\}+uu\{ \{ c,a\},b\}]=0.$$

The theorem is proved.

\subsection{Generalized Poisson algebras.}
The variety of unital {\it generalized Poisson algebras} (see, for example, \cite{kaysh14}) is
 defined by the identities
$$xy=yx, (xy)z=x(yz), \{xy,z\}=\{x,z\}y+x\{y,z\}+D(z)xy,D(x)=\{x,1\},$$
$$\{x,y\}=-\{y,x\}, \{\{x,y\},z\}+\{\{y,z\},x\}+\{\{z,x\},y\}=0.$$

    \begin{Th}\label{thm_genpoisson}Let $(A;\cdot,\{,\})$ be a generalized Poisson algebra.
    Then  $(A;[\{,\}, \cdot])$ is an associative-commutative algebra, and
    $(A;[\cdot, \{,\}])$ is a Lie  algebra.
    \end{Th}

{\bf Proof.}
In the first case, we have
$$a*b= \{u,ab\} - \{u,a\}b-a\{u,b\}=-D(u)ab=b*a$$
and
$$(a*b)*c=D(u)^2abc=a*(b*c).$$

In the second case,
$$a*b= u\{a,b\}-\{ua,b\}-\{a,ub\}=$$
$$- a\{u,b\}- b\{a,u\}-u\{a,b\} - D(b)ua+D(a)ub =-b*a.$$

Here we use the proof of Theorem \ref{ThPoisson}.
It is easy to see that
\begingroup\small\begin{eqnarray*}
\circlearrowright_{a,b,c} [&& (a*b)*c]= \circlearrowright_{a,b,c} [(\{b,u\}a+\{b,a\}u+
\{u,a\}b-D(b)ua+D(a)ub)*c]=\\
\circlearrowright_{a,b,c} [&&\{c,u\}\{b,u\}a+\{c,b\}\{b,u\}a+\{c,\{b,u\}\}au-
\{b,u\}D(c)au+\{u,a\}\{b,u\}c+\\
&&\{u,\{b,u\}\}ac-\{b,u\}D(u)ac-\{b,u\}D(c)au+\{b,u\}D(a)cu+D(\{b,u\})acu+\\
&&\{c,u\}\{b,a\}u+\{c,u\}\{b,a\}u+\{c,\{b,a\}\}u^2-\{b,a\}D(c)u^2+\{u,\{b,a\}\}cu-\\
&& \{b,a\}D(u)cu-\{b,a\}D(c)u^2+\{b,a\}D(u)cu+D(\{b,a\})cu^2+\{c,u\}\{u,a\}b+\\
&& \{c,b\}\{u,a\}u+\{c,\{u,a\}\}bu-\{u,a\}D(c)bu+\{u,\{u,a\}\}bc+\{u,b\}\{u,a\}c-\\
&& \{u,a\}D(u)bc-\{u,a\}D(c)bu+\{u,a\}D(b)cu+D(\{u,a\})bcu-\{c,u\}D(b)au-\\
&& \{c,D(b)\}au^2-\{c,u\}D(b)ua-\{c,a\}D(b)u^2+2D(b)D(c)au^2-\{u,a\}D(b)cu-\\
&& \{u,D(b)\}acu+2D(b)D(u)acu+D(b)D(c)au^2 -D(D(b))acu^2-D(b)D(u)acu-\\
&& D(a)D(b)acu^2+\{c,u\}D(a)bu+\{c,D(a)\}bu^2+\{c,u\}D(a)bu+\{c,b\}D(a)u^2-\\
&& 2D(a)D(c)bu^2+\{u,D(a)\}bcu+\{u,b\}D(a)cu-2D(a)D(u)bcu-D(a)D(c)bu^2+\\
&& D(D(a))bcu^2+D(a)D(u)bcu+D(a)D(b)cu^2].\
\end{eqnarray*}\endgroup{}

By the proof of Theorem \ref{ThPoisson}, we can conclude that the sum of all elements without $D$ is zero.
Now, we have
\begingroup\small\begin{eqnarray*}
\circlearrowright_{a,b,c} [&& (a*b)*c]= \\
\circlearrowright_{a,b,c} [&&
-\{b,u\}D(c)au-\{b,u\}D(u)ac-\{b,u\}D(c)au+\{b,u\}D(a)cu+D(\{b,u\})acu-\\
&& \{b,a\}D(c)u^2-\{b,a\}D(u)cu-\{b,a\}D(c)u^2+\{b,a\}D(u)cu+D(\{b,a\})cu^2-\\
&& \{u,a\}D(c)bu-\{u,a\}D(u)bc-\{u,a\}D(c)bu+\{u,a\}D(b)cu+D(\{u,a\})bcu-\\
&& \{c,u\}D(b)au-\{c,D(b)\}au^2-\{c,u\}D(b)ua-\{c,a\}D(b)u^2+2D(b)D(c)au^2-\\
&& \{u,D(b)\}acu+2D(b)D(u)acu+D(b)D(c)au^2 -D(D(b))acu^2-D(b)D(u)acu-\\
&& \{u,a\}D(b)cu-D(a)D(b)acu^2+\{c,u\}D(a)bu+\{c,D(a)\}bu^2+\{c,u\}D(a)bu+\\
&& \{c,b\}D(a)u^2-2D(a)D(c)bu^2+\{u,D(a)\}bcu+\{u,b\}D(a)cu-2D(a)D(u)bcu-\\
&& D(a)D(c)bu^2+D(D(a))bcu^2+D(a)D(u)bcu+D(a)D(b)cu^2].\\
\end{eqnarray*}\endgroup{}

It is easy to see that
\begin{eqnarray*}
\circlearrowright_{a,b,c} && [D(D(a))bcu^2-D(D(b))acu^2]=0,\\
\circlearrowright_{a,b,c} && [D(\{b,u\})acu+D(\{u,a\})bcu]=0,\\
\circlearrowright_{a,b,c} && [\{u,D(a)\}bcu-\{u,D(b)\}acu]=0,\\
\circlearrowright_{a,b,c} && [D(\{b,a\})cu^2-\{c,D(b)\}au^2+\{c,D(a)\}bu^2]=0,\\
\circlearrowright_{a,b,c} && [-2\{b,a\}D(c)u^2-\{c,a\}D(b)u^2+\{c,b\}D(a)u^2]=0,\\
\circlearrowright_{a,b,c} && [\{b,u\}D(c)au+\{u,a\}D(c)bu+\{c,u\}D(b)au+\{u,c\}D(a)bu]=0.
\end{eqnarray*}

Obviously, $\circlearrowright_{a,b,c} [(a*b)*c]=0$
and  $[\cdot, \{,\}]$ is a Lie  algebra.
The theorem is proved.

\

\subsection{Novikov-Poisson algebras.}
The variety of left {\it Novikov-Poisson} algebras is defined by the identities
$$xy=yx, (xy)z=x(yz),$$
$$x \circ (y \circ z)= y \circ (x \circ z), (x, y, z)_{\circ}  =  (x,z, y)_{\circ},$$
$$x\circ (y  z) = (x  \circ  y)  z,
(x   y) \circ  z - x   (y \circ  z)= (x  z) \circ y - x   (z \circ y).$$

    \begin{Th}Let $(A;\cdot,\circ)$  be a left Novikov-Poisson algebra.
    Then  $(A;[\cdot, \circ])$ is a left Novikov algebra and
          $(A;[\circ, \cdot])$ is an associative-commutative  algebra.
    \end{Th}

{\bf Proof.}
Firstly, we have
$$a * b = u (a \circ b)-(ua) \circ b-a \circ (ub)= -(ua)\circ b.$$
 Hence,
$$ a*(b*c)= (ua)\circ ((ub)\circ c)= (ub)\circ ((ua)\circ c)= b*(a*c),$$
and
$$(a,b,c)_*=(a*b)*c- a*(b*c)=
(u((ua)\circ b))\circ c - (ua)\circ ((ub)\circ c)=$$
$$((ua)\circ (ub))\circ c - (ua)\circ ((ub)\circ c) =
(ua, ub, c)_{\circ}=(ua, c, ub)_{\circ}=$$
$$=((ua)\circ c)\circ(ub)-(ua)\circ(c\circ(ub))=$$
$$u(ua,c,b)_{\circ}=u(ua,b,c)_{\circ}=(a*c)*b-a*(c*b)=(a,c,b)_*.$$
Secondly,
$a * b = u \circ (ab)-(u\circ a) b-a(u\circ b)=-u\circ(ab)=b*a.$
Therefore,
$$(a *  b) * c= u \circ ( (u \circ (ab)) c)=u \circ  (u \circ (abc))=u \circ (a (u \circ (bc)))= a*(b*c).$$
The theorem is proved.

\

Similarly, the variety of right Novikov-Poisson algebras may be defined  (see, for example, \cite{xu97}).
 It is easy to prove the following theorem:

    \begin{Th}Let $(A;\cdot,\circ)$  be a right Novikov-Poisson algebra.
    Then  $(A;[\cdot, \circ])$ is a right Novikov algebra and
          $(A;[\circ, \cdot])$ is a commutative  algebra.
    \end{Th}

\section{The Kantor square of algebras of special type.}

Here we study some special cases of the Kantor square.
For an algebra $A=(A;\cdot)$, we denote its Kantor square $(A;[\cdot,\cdot])$ by $(A,*)$.
We denote the Kantor square for a fixed element $u$ by $(A,*_u)$.
We consider the relations between the ideals in $A$ and $(A,*)$,
the relations between an associative algebra $A$ with a polynomial identity and its Kantor square.
Moreover, the relations between the nilpotency and right nilpotency in $A$ and $(A,*)$ are  investigated.

\subsection{Ideals in the Kantor product.}

\begin{Th}\label{thm_ideal}
Let $I$ be an ideal of  $A$. Then $(I,*)$ is an ideal of $(A,*),$
but the converse statement  is not true in general.
\end{Th}

{\bf Proof.}
It is easy to see that
if $i \in I$ and $a \in A$ then
$$i*a=u(ia)-(ui)a-i(ua) \in I \mbox{ and }a*i=u(ai)-(ua)i-a(ui) \in I.$$
It  follows that $I$ is an ideal of $(A,*).$

Conversely,
we can consider the trivial case
where $A$ has zero Kantor square (for example, a Lie or Leibniz algebra)
and every subspace of $A$ is an ideal of $(A,*).$
For the non-trivial case (nonzero Kantor product),
we can consider the following associative algebra: $A_1 \oplus A_2$
 is the direct sum of the matrix algebras of order $2$.
 Here, if $e_i$ is the unit of $A_i$ then
 the subspace generated by $A_1$ and $e_2$ is an ideal of $(A,*_{e_1}),$
but is not an ideal of $A_1 \oplus A_2.$
The theorem is proved.

\subsection{Associative algebras with polynomial identity.}
Given a polynomial $f$ in $n$ variables,
we define $f_*(x_1, \ldots, x_n)$
as the value of $f$ in $(A,*)$, where $x_1, \ldots, x_n$ are some elements in $ A.$

\begin{Th}\label{ass_id1}
Let $(A;\cdot)$ be an associative algebra that satisfies the polynomial identity $f(x_1,\ldots, x_n)$.
 Then there exists an identity $g$ such that
 $A$ and $(A,*)$ satisfy $g$.
\end{Th}

{\bf Proof.}
Since every algebra satisfying a polynomial identity also satisfies a nonzero multilinear identity of at most the same degree (obtained by the standard linearization process, see~\cite{rings}), we may assume that $f(x_1, \ldots, x_n)$ is multilinear.
Then $A$ satisfies the identity $g(x_1, \ldots, x_n,z)=f(x_1, \ldots, x_n)z.$
By the definition of the Kantor product (see Lemma~\ref{thass}),
the multiplication in the algebra $(A,*)$ is defined by $x*y=-xuy.$
Now, we can calculate the element $g_*(x_1, \ldots, x_n,z)$ in $(A,*)$.
Since every monomial of $f$ has degree $n$, it is $(-1)^{n}f(x_1u, \ldots, x_nu)z$, which is zero in $A$.
It follows that $(A,*)$ satisfies the identity $g.$
The theorem is proved.

One of the most popular identities in associative algebras is the standard polynomial identity of degree $n$:
$$s_{n}(x_1, \ldots, x_n)= \sum_{\sigma \in S_{n}} (-1)^{\sigma} x_{\sigma(1)} \ldots  x_{\sigma(n)}.$$

\begin{Th}\label{ass_id2}
Let $(A;\cdot)$ be an associative algebra that
 satisfies $s_n$. Then $(A,*)$ satisfies $s_{n+1}.$

\end{Th}

{\bf Proof.}
It is easy to see that
 the standard polynomial of degree $n+1$ may be written as
$$ s_{n+1}(x_1, \ldots, x_{n+1})=
\sum_{i=1}^n (\sum_{ \begin{array}{c} \small\sigma \in S_{n+1}, \\ \sigma(n+1)=i
\end{array}} (-1)^{\sigma} x_{\sigma(1)} \ldots  x_{\sigma(n)} x_{\sigma(n+1)=i})=$$
$$=\sum_{i=1}^n (  \epsilon_{i} s_{n}(x_1, \ldots, \hat{x_i}, \ldots, x_{n+1}) x_{i}),
\epsilon_{i} = \pm 1.$$
Now, by the proof of Theorem \ref{ass_id1},  $(A,*)$ satisfies
 the standard polynomial identity of degree $n+1.$
The theorem is proved.

\subsection{Nilpotent algebras.}
For nilpotent algebras, we can prove the following lemma.

\begin{Lem}\label{nilpotent}
Let $(A;\cdot)$ be a nilpotent algebra of nilpotency index $n$.
Then $(A,*)$ is a nilpotent algebra of nilpotency index $\leq [n/2]+1.$
\end{Lem}

{\bf Proof.}
Obviously, every product of the form  $x_1*x_2*\ldots *x_t$ (with some order of brackets)
is a sum of
 multiplications of the form $y_1y_2\cdots y_{2t-1}$ (with some order of brackets).
 Now, it is easy to see that $(A,*)$ is nilpotent and its index of nilpotency is $\leq [n/2]+1.$
The lemma is proved.

\subsection{Right-nilpotent algebras.}
An algebra $A$ is called right-nilpotent (or left-nilpotent) of nilpotency index $n$
 if it satisfies the identity
$$(\ldots (x_1x_2)\ldots)x_n=0 \ \ \ (\mbox{ or }   x_1( \ldots (x_{n-1}x_n)\ldots )=0).$$
Curiously, an analogue of Lemma~\ref{nilpotent} is not true for right-nilpotent algebras.

\begin{Th}\label{rightnilpotent}
There exists a right nilpotent algebra
 $(A;\cdot)$  such that $(A,*_u)$ is not right nilpotent, but  $(A,*_u)$  is solvable.
\end{Th}

{\bf Proof.} An algebra $A$ is \textit{right alternative} if the following identity holds in $A$:
$(x,y,y) = 0.$
It is an interesting fact that, in contrast to the algebras of many well-studied classes
 (Jordan, alternative, Lie, and so on), a right nilpotent right alternative algebra need not be nilpotent.
The corresponding example of a five-dimensional right nilpotent but not nilpotent algebra belongs to
  Dorofeev \cite{Dorofeev}.
Its basis is $\{a, b, c , d , e \},$  and the multiplication is given
 by (zero products of basis vectors are omitted)
$$ab = -ba = ae = -ea = db = -bd = -c, ac = d, bc = e.$$
It is easy to see that
$c *_ab = a(cb)-(ac)b-c(ab)=c.$
Obviously, $c=( \ldots (c*_a b)*_a \ldots )*_ab \neq 0,$
and $(A,*_a)$ is not right-nilpotent.
It is easy to see that $A^2 \subseteq \langle c,d,e \rangle$ and $A*_aA  \subseteq \langle c,d,e \rangle$,
 but $(A*_aA)*_a(A*_aA)=0,$
 and $(A,*_a)$ is solvable.

The theorem is proved.

\subsection{Derivations.}
Remember that a linear mapping $D$ of an algebra $A$ is called a derivation
if it satisfies the relation $D(xy)=D(x)y+xD(y).$
Following~\cite{Kantor90}, an element $a$ of an algebra $A$ is called a {\it Jacobi} element if
it satisfies the relation $a(xy)=(ax)y+x(ay).$
All Jacobi elements of $A$ form a vector space, which is  called the {\it Jacobi} space of $A.$

\begin{Lem}\label{derivation}
Let $D$ be a derivation of both $A$ and $(A,*_u).$
Then $D(u)$ is a Jacobi element of $A$. In particular,

$1)$ if $A$ has zero Jacobi space, then $D(u)=0$; hence, if $A$ has zero Jacobi space and $D$ is a derivation of $(A,*_u)$ for every $u\in A$, then $D=0;$

$2)$ if $D$ is invertible and $D$ is a derivation of $(A,*_u)$ for every $u\in A$, then $A$ is a left Leibniz algebra and $(A,*_u)$ is a zero algebra for every $u\in A$.
In particular, if $A$ is a finite-dimensional algebra over a field of characteristic zero, then $A$ is nilpotent.

\end{Lem}

{\bf Proof.}
By simple calculations, from
$D(x*y)=D(x)*y+x*D(y)$ and $D(xy)=D(x)y+xD(y),$
we have
$$D(u)(xy)=(D(u)x)y+x(D(u)y),$$
i.e. $D(u)$ is a Jacobi element of $A$. This gives $1)$.

$2).$ If $D$ is a derivation of $(A,*_u)$ for every $u\in A$, then every element of $D(A)=A$ is a Jacobi element, i.e. $A$ is a left Leibniz algebra.
By Lemma~\ref{leibniz}, $(A,*_u)$ is a zero algebra for every $u\in A$.
In~\cite{Fialc12} it was proved that a finite-dimensional Leibniz algebra
over a field of characteristic zero admitting an invertible derivation is nilpotent.

Note that in $2)$ the condition that $D$ is a derivation of $(A,*_u)$ for every $u\in A$ cannot be omitted. Indeed, let $A$ be the algebra with a basis $e_1,e_2,e_3$ and the non-zero products $e_1e_1=e_2$, $e_2e_1=e_3$. Then $D(e_i)=ie_i$ is an invertible derivation of $A$ and of the zero algebra $(A,*_{e_3})$, but $A$ is not a left Leibniz algebra, since $e_1(e_1e_1)=0\neq e_3=(e_1e_1)e_1+e_1(e_1e_1)$.
The lemma is proved.

\subsection{Automorphisms.}
Remember that an invertible linear mapping $\phi$ of an algebra is called an automorphism
if it satisfies the relation $\phi(xy)=\phi(x)\phi(y).$

\begin{Lem}\label{automorphism}
Let $\phi$ be an automorphism of both $A$ and $(A,*_u).$
Then $u-\phi(u)$ is a Jacobi element of $A$. In particular, if $A$ has zero Jacobi space, then $\phi(u)=u$; hence, if $A$ has zero Jacobi space and $\phi$ is an automorphism of $(A,*_u)$ for every $u\in A$, then $\phi$ is the identity mapping.
\end{Lem}

{\bf Proof.}
By simple calculations, from
$\phi(x*y)=\phi(x)*\phi(y)$ and $\phi(xy)=\phi(x)\phi(y),$
we have
$$(u-\phi(u))(\phi(x)\phi(y))=((u-\phi(u))\phi(x))\phi(y)+\phi(x)((u-\phi(u))\phi(y)).$$
Since $\phi$ is surjective, $u-\phi(u)$ is a Jacobi element of $A$.
The lemma is proved.

\subsection{Isomorphic Kantor squares.}
Here we talk about the situation where  algebra $A$ and its Kantor square are isomorphic.

\begin{Th}\label{Isomor}
Let $A$ be a unital associative algebra and $u\in A$.
Then $A$ is isomorphic to $(A,*_u)$
if and only if $u$ is invertible. In particular, if $A$ is a skew field, then $A$ is isomorphic to $(A,*_u)$ for every nonzero $u\in A$.

\end{Th}

{\bf Proof.}
Recall that $x*_uy=-xuy$ (see the proof of Theorem~\ref{ass_id1}).
If $u$ is invertible, then the map $f_u(a)=-au^{-1}$ is an isomorphism between the algebras $A$ and $(A,*_u)$, since
$$f_u(a)*_uf_u(b)=-(-au^{-1})u(-bu^{-1})=-abu^{-1}=f_u(ab).$$
Conversely, let $f:A\to(A,*_u)$ be an isomorphism and $e=f(1)$. Then $e$ is the unit of $(A,*_u)$, i.e. $-euy=y=-yue$ for all $y\in A$. Taking $y=1$, we obtain $(-e)u=1=u(-e)$, i.e. $u$ is invertible.
In particular, the matrix algebra $M_n(F)$ $(n\geq2)$ is isomorphic to $(M_n(F),*_1)$, but it is not a skew field.
The theorem is proved.

\

\appendix
\section*{Appendix A. Differences from the published version}
\addcontentsline{toc}{section}{Appendix A. Differences from the published version}

\noindent\textbf{Lemma~\ref{automorphism}:} the statement has been changed (in general only $u-\phi(u)$ is a Jacobi element) and the proof corrected.

\noindent\textbf{Theorem~\ref{Isomor}:} the statement has been changed and the proof corrected: the original statement is false, since $M_2(F)\cong(M_2(F),*_1)$.

\end{document}